\documentclass[11pt]{article}

\usepackage{amssymb,amsmath,amsthm,rotating}

\newcommand{\n}{\noindent}
\newcommand{\bb}[1]{\mathbb{#1}}
\newcommand{\cl}[1]{\mathcal{#1}}
\newcommand{\mf}[1]{\mathfrak{#1}}
\newcommand{\lngarr}{\longrightarrow}
\newcommand{\ds}{\displaystyle}

\theoremstyle{plain}
\newtheorem{pro}{Proposition}[section]
\newtheorem{lem}{Lemma}[section]
\newtheorem*{thmb}{Theorem}

\begin{document}

\title{Ideals in Toeplitz Algebras\thanks{This paper was written during a visit at the Fields
Institute in Toronto while on research leave from Texas A\&M University.}}

\author{Ronald G.\ Douglas}

\date{}
\maketitle

\begin{abstract}
 We determine the ideal structure of the Toeplitz $C^*$-algebra 
on the bidisk
\end{abstract}

\setcounter{section}{-1}

\section{Introduction}\label{sec0}

\indent

A large part of doing research in mathematics is asking the right question. 
Posing a timely question can trigger thought and provoke insights, even when the 
matter has no good resolution. That is what happened in the case at hand. After 
devoting much time to the study of Hilbert modules and quotient Hilbert modules 
in particular, I had occasion to wonder if a $C^*$-algebra perspective on this 
topic might not be useful. More precisely, I asked myself whether there was any 
relationship between the $C^*$-algebras for  a Hilbert module and its quotient 
module  generated by the operators obtained from module multiplication. 
The answer turned out to be negative but the techniques and approach used to 
reach that conclusion, namely an analysis of the ideal structure of the Toeplitz 
$C^*$-algebras, seemed worth reporting on. Although this general topic has been 
considered 
earlier (cf.\ \cite{upmeier}), these studies focused almost exclusively on index 
theory 
which will not concern us here.

Let $\Omega$ be a bounded domain in $\bb C^m$ and $A(\Omega)$ be the function 
algebra obtained from the closure in the supremum norm on $\Omega$ of all 
functions holomorphic on some neighborhood of the closure of $\Omega$. In this 
paper we will focus on the case of the unit ball $\bb B^m$ and the polydisk $\bb 
D^m$ for $m=2$.
\vfill

\n {\bf 2000 Mathematical Subject Classification.} 46L06, 47L05, 47L15, 47L20, 
47L80.

\n {\bf Key words and phrases.}  Toeplitz operators, Toeplitz $C^*$-algebras, 
Hardy 
space on bidisk, ideals in $C^*$-algebras.
\newpage

A Hilbert module $\cl M$ over $A(\Omega)$ is a Hilbert space together with a 
unital module map $A(\Omega\times \cl M\to \cl M$ which is continuous. Here 
we confine attention to contractive Hilbert modules; that is, a module 
$\cl M$ for which the inequality $\|\varphi f\|_{\cl M} \le 
\|\varphi\|_{A(\Omega)} \|f\|_{\cl M}$ holds for $\varphi$ in $A(\Omega)$ and 
$f$ in $\cl M$. Actually, we will be concerned only with the Hardy modules 
$H^2(\bb D^m)$ and $H^2(\bb B^m)$ on $\bb D^m$ and $\bb B^m$, respectively, and 
their quotient modules. Recall that the Hilbert space for this module can be 
obtained as the closure of $A(\bb D^m)$ in the space $L^2((\partial \bb D)^m)$ 
and $L^2(\bb B^m)$ relative to Lebesgue measure on the distinguished boundary 
$(\partial\bb D)^m$ of $\bb D^m$ and the boundary $\partial \bb D^m$,
respectively.

In \cite{douglas2} quotient modules obtained from a hypersurface $\mathfrak{Z}$ 
are studied 
and a characterization is obtained of them based on the intrinsic geometry of 
how $\mathfrak{Z}$ sets in $\Omega$ and the spectral sheaf for the ambient 
Hilbert module. We will not provide any details on those results here except to 
point out that their precision might suggest that an intimate relationship might 
also exist for the corresponding $C^*$-algebras. In particular, let 
$\mathfrak{C}^*(\cl 
M)$ denote the $C^*$-algebra generated by the operators on $\cl M$ defined by 
module  multiplication. The question which spawned this study is whether there 
is a relationship between $\mathfrak{C}^*(\cl R)$ and $\mathfrak{C}^*(\cl R/\cl 
R_0)$ for $\cl R$ the Hardy module and $\cl R_0$ the submodule of functions in 
$\cl R$ that 
vanish on $\mf{Z}$. Another way to phrase the question is whether the map 
defined by $M_\psi\to N_\psi$ for $\psi$ in $A(\Omega)$ extends to the 
$C^*$-algebras $\mf C^*(R)$ and $\mf C^*(\cl R/\cl R_0)$, where $M_\psi$ and 
$N_\psi$ are the operators defined by module 
multiplication  on $\cl R$ and $\cl R/\cl R_0$, respectively.

\section{Quotient Modules Determined by Hypersurfaces}\label{sec1}

\indent

One way to approach this question is by first determining the ideal $\mf I$ in 
$\mf C^*(\cl R)$ generated by $\{M_\varphi\mid \varphi|_{\mf Z} \equiv  0\}$, 
which leads 
to 
the topic of this paper:\ What are the ideals in $\mf C^*(\cl R)$? More 
particularly, 
determine the ideal generated by a collection of operators and calculate the 
quotient $C^*$-algebra. For the case of the ball,  $\mf C^*(H^2(\bb 
B^m))$ contains the ideal $\mf K(H^2(\bb B^m))$ of compact operators on $H^2(\bb 
B^m)$ and the quotient is commutative and isomorphic to $C(\partial \bb B^m)$. 
Thus the ideal structure of $\mf C^*(H^2(\bb B^m))$ is especially simple. 
The ideals in $\mf C^*(H^2(\bb B^m))$ are  $\mf K(H^2(\bb B^m))$ and those 
corresponding to  closed subsets of $\partial\bb B^m$. As a result, there can be 
no map from $\mf C^*(H^2(\bb B^m))$ to the $C^*$-algebra determined by a 
quotient module since the latter algebras are not commutative. The situation for 
$\mf C^*(H^2(\bb D^m))$, however, is more interesting. Rather than tackling this 
problem in full generality at the outset, however, we begin by analyzing some 
examples.

Let us consider three hypersurfaces in the closure of $\bb D^2$; the quotient 
modules they determine; and the $\mf C^*$-algebra maps: 
\begin{align*}
\mf Z_1 &= \{(z_1,z_2)\mid z_1=0, |z_2| \le 1\},\\
\mf Z_2 &= \{(z_1,z_2)\mid z_1 = 1, |z_2| \le 1\} \text{ and}\\
\mf Z_3 &= \{(z_1,z_2)\mid z_1=z_2, |z_1| \le 1, |z_2|\le 1\}.
\end{align*}
 Using the results in \cite{douglas2} or by explicit calculation, one can show 
the quotient $C^*$-algebras are:
\begin{align*}
\mf Q _1 &= H^2(\bb D^2)/\{f\in H^2(\bb D^2)\mid f|_{\mf Z_1} \equiv 0\} \cong 
H^2(\bb D),\\
\mf Q_2 &= H^2(\bb D^2)/\{f\in H^2(\bb D^2)\mid f|_{\mf Z_2} \equiv 0\} = (0), 
\text{ and }\\
\mf Q_3 &= H^2(\bb D^2)/\{f\in  H^2(\bb D^2)\mid f|_{\mf Z_3} \equiv 0\} \cong 
B^2(\bb D).
\end{align*}
(Here, $B^2(\bb D)$ is the Hilbert module over the disk algebra $A(\bb D)$ 
determined by the Bergman space on $\bb D$.)

To make sense of these equivalences, we must provide a few more details. In the 
first case, the associated module multiplication on $H^2(\bb D)$ has $z_1$ 
acting as the zero operator while $z_2$ is just multiplication by $z$ on 
$H^2(\bb D)$. In the second case, we interpret the condition $f|_{\mf Z_2}\equiv 
0$ 
to select the set of functions in $A(\bb D^2)$ that vanish on $\mf Z_2$. Hence, 
the closure in $H^2(\bb D^2)$ of such functions is all of $H^2(\bb D^2)$ and 
thus the quotient module is (0). Finally, the quotient $\mf  Q_3$ is 
isometrically isomorphic to the Bergman space with both multiplication by $z_1$ 
and $z_2$ corresponding to multiplication by $z$. (This equivalence was first 
observed by Rudin \cite{rudin}.) To examine whether the maps from $M_\varphi$ to 
$N_\varphi$ for $\varphi$ in $A(\bb D^2)$ extend to the $C^*$-algebras, we want 
to consider the ideals $\mf I_1$, $\mf I_2$ and $\mf I_3$ in $\mf C^*(H^2(\bb 
D^2))$ generated by the sets $\{M_\varphi\mid \varphi|_{\mf Z_i}\equiv 0\}$ for 
$i=1,2,3$. We state and prove our general theorem characterizing ideals after 
working out these examples.

\section{Ideals in Toeplitz $\pmb{C^*}$-Algebras---Three Examples}\label{sec2}

\indent

We begin by recalling something about the nature of the various Toeplitz 
$C^*$-algebras. For $H^2(\bb D^2)$ we now denote it by $\mf T (H^2(\bb D^2))$. 
This 
bidisk Toeplitz $C^*$-algebra was introduced by R.\ Howe and the author in 
\cite{douglas1} 
where it was used to study the invertibility and Fredholmness of Toeplitz 
operators in the quarter-plane. We'll recall the techniques developed in that 
paper a little later. Next,  the standard Toeplitz $C^*$-algebra 
$\mf T(H^2(\bb D))$ for the disk is obtained by considering module 
multiplication on the 
quotient module $H^2(\bb D^2)/[z_1]$, where $[z_1]$ denotes the closed submodule 
generated by $z_1$. In the second case, asking about the algebra for the 
quotient module is vacuous since $H^2(\bb D^2)/[z_1-1] = (0)$. Finally, the 
$C^*$-algebra defined by module multiplication on the Bergman module is known to 
be $\mf T(H^2(\bb D))$, the same as the Toeplitz $C^*$-algebra for $H^2(\bb D)$. 
This result follows from an explicit calculation but also from the BDF 
characterization of essentially normal operators \cite{brown}.

To study the maps on $\mf T(H^2(\bb D^2))$ defined by $M_\varphi\to 
N_{\varphi|\mf Z_i}$ 
for $\varphi$ in $A(\bb D^2)$, we begin by considering the ideal $\mf I_i$ in 
$\mf T(H^2(\bb D^2))$ generated by the Toeplitz operators whose symbols vanish 
on 
$\mf Z_i$ for $i=1,2,3$. Recall that corresponding to the identification of 
$H^2(\bb D^2)$ with the Hilbert space tensor product $H^2(\bb D) \otimes H^2(\bb 
D)$, we have $\mf T(H^2(\bb D^2))\cong \mf T(H^2(\bb D))$ $\otimes~\mf T(H^2(\bb 
D))$, 
where the latter is the spatial tensor product. Moreover, we use the standard 
notation of $T_\varphi$ for the Toeplitz operator defined by module 
multiplication by $\varphi$. Finally, we let $z$ denote the variable on $\bb D$ 
and $z_1$ and $z_2$ denote the variables on $\bb D^2$.

\subsection{Example 1}

\indent

For $\mf I_1$, we note that $T_{z_1} = T_z\otimes I$ and $T^*_{z_1} T_{z_1} = 
(T_z \otimes I)^* (T_z\otimes I) = I\otimes I = I$ and hence the identity 
operator is in $\mf I_1$ which means that $\mf I_1 =  \mf T(H^2(\bb D^2))$. 
Therefore, the map from module multipliers on $H^2(\bb D^2)$ to their 
restrictions acting on the quotient $H^2(\bb D^2)/[z_1]$ does not extend.

\subsection{Example 2}

\indent

Before tackling the next example we need to recall the basic results from 
\cite{douglas1} 
in which the short exact sequence for the Toeplitz $C^*$-algebra on the disk
\[ 
0\to \mf K(H^2(\bb D)) \lngarr  \mf T(H^2(\bb D)) 
\overset{\ds\sigma}{\longrightarrow} 
C(\partial \bb D) \lngarr 0
\]
is tensored with itself. (In \cite{berger} it was pointed out that while the 
preservation 
of exactness for tensor products is not automatic, it is valid in this case and, 
moreover, all the tensor products are unique.) We will use $\mf K$ and $\mf T$ 
to 
denote the $C^*$-algebras $\mf K(H^2(\bb D))$ and $\mf T(H^2(\bb D))$, 
respectively.

\begin{pro}\label{pro2.1}
The following diagram is exact:
\newpage

\[
\begin{matrix}
&&0&&0&&0\\
&&\downarrow&&\downarrow&&\downarrow\\
0&\lngarr&\mf K\otimes \mf K&\lngarr&\mf K\otimes \mf T&\lngarr&\mf K\otimes 
C(\partial D)&\lngarr&0\\
&&\downarrow&\raisebox{3ex}{\begin{rotate}{-45}{$\overset{\ds i}{\hbox to 
35pt{\rightarrowfill}}$}\end{rotate}}\hfill&\downarrow&&\downarrow\\
&&&&&&\\
0&\lngarr&\mf T\otimes\mf K&\lngarr&\mf T\otimes \mf T&\lngarr&\mf T\otimes 
C(\partial\bb 
D)&\lngarr&0\\
&&\downarrow&&\downarrow&\raisebox{.5ex}{\begin{rotate}{-45}{$\overset{\ds 
\Sigma}{\hbox to 35pt{\rightarrowfill}}$}\end{rotate}}~~~~\hfill&\downarrow\\
&&&&&&\\
0&\lngarr&C(\partial D)\otimes\mf K&\lngarr&C(\partial \bb D)\otimes \mf T&\lngarr 
&C(\partial\bb D)\otimes C(\partial\bb D))&\lngarr&0\\
&&\downarrow&&\downarrow&&\downarrow\\
&&0&&0&&0
\end{matrix}
\]
\end{pro}

We observe that spatially we have $\mf K\otimes \mf K \cong \mf K(H^2(\bb D^2))$ 
and $\mf T \otimes \mf T\cong  \mf T(H^2(\bb D^2))$. Moreover, one has 
$C(\partial \bb D) 
\otimes C(\partial\bb D)\cong C((\partial \bb D)^2)$ and $\mf K\otimes 
C(\partial\bb D)\cong C(\partial \bb D)\otimes \mf K \cong C(\partial\bb D,\mf 
K)$ and $\mf T\otimes C(\partial\bb D)\cong C(\partial\bb D) \otimes  \mf 
T\cong 
C(\partial \bb D,\mf T)$, where $C(\partial\bb D,\mf K)$ and $C(\partial \bb D, 
\mf T)$ 
are the $C^*$-algebras of continuous functions on $\partial\bb D$ taking values 
in $\mf K$ and $\mf T$, respectively. Finally, we record for later use:

\begin{pro}\label{pro2.2}
If $i$ and $\Sigma$ are defined by the diagram, then
\begin{align*}
\text{Range } i &= \text{Kernel}(\mf T \otimes \mf T \to (\mf T \otimes 
C(\partial\bb D) 
\oplus C(\partial\bb D)\otimes \mf K) \text{ and }\\
\text{Kernel } \Sigma &= \text{Range}((\mf T \otimes \mf K \oplus \mf K \otimes 
 \mf T) 
\lngarr  \mf T \otimes \mf T).
\end{align*}
\end{pro}

The proof involves a standard diagram chase.

Now $\mf I_2$ is the ideal in $\mf T(H^2(\bb D^2))$ generated by the operator 
$T_{z_1-1}$, since the set of functions of the form $(z_1-1) \varphi(z_1,z_2)$ 
is uniformly dense in the ideal of functions in $A(\bb D^2)$ that vanish on 
$\mf Z_2$. Calculating, we have
\[
\Sigma(T_{z-1}\otimes I) = (z-1)\otimes 1 = z_1 -1
\]
and therefore
\[
\Sigma(\mf I_2) = \{f\in C((\partial \bb D)^2)\mid f|_{\mf Z_2}\equiv 0\}.
\]
Further, the commutator $[T^*_{z-1}\otimes I,T_{z-1}\otimes I] = P_0\otimes I$ 
is in $(\mf K\otimes \mf T)\cap \mf I_2$, where $P_0$ is the projection onto the 
constant functions in $H^2(\bb D)$. Since $\mf I_2$ is an ideal, we see that 
$P_0\otimes \mf T$ is contained in $\mf I_2$ and hence $\mf K\otimes  \mf T 
\subset \mf 
I_2$. Now appealing to the exact sequence $0\to \mf K\otimes \mf T \to  \mf T 
\otimes \mf T \overset{\ds \sigma\otimes 1}{\hbox to 35pt{\rightarrowfill}} 
C(\partial \bb D, 
\mf T)\to 0$, we see since $\mf K\otimes  \mf T \subset \mf I_2$ that $\mf I_2$ 
is 
determined by $(\sigma\otimes 1)(\mf I_2)$. Therefore, $\mf I_2 = \{X\in  
\mf T\otimes \mf T|(\sigma\otimes 1)(X)|_{z=1} \equiv 0\}$.

Now consider the following diagram obtained by taking quotients:
\[
\begin{matrix}
\mf T\otimes \mf T&\lngarr&C(\partial \bb D, \mf T)\hfill\\
\Big\downarrow&&\Big\downarrow\hfill\\
\mf T \otimes  \mf T/[T_{z_1-1}]&\lngarr&C(\partial\bb D, \mf T)/[z_1-1] \cong  
\mf T(H^2(\bb D)).\hfill
\end{matrix}
\]
Using the bottom line, we see that  the map $T_\varphi\to T_{\varphi|\mf Z_1}$ 
defined for $\varphi$ in $A(\bb 
D^2)$ extends to a $C^*$-homomorphism from $\mf T(H^2(\bb D^2))$ to $\mf 
T(H^2(\bb 
D))$. Thus, although the underlying Hilbert modules don't match up in this case, 
the $C^*$-algebras generated by the module multipliers do.

\subsection{Example 3}

\indent

Now let us consider $\mf I_3$ which is the ideal in $\mf T(H^2(\bb D^2))$ 
generated 
by $T_{z_1-z_2}$ or $T_z\otimes I - I\otimes T_z$ in $\mf T\otimes  \mf T$. We 
proceed 
as in the last case by calculating the commutator
\[
[(T_z\otimes I-I\otimes T_z)^*, T_z\otimes I-I\otimes T_z] = P_0 \otimes I - I 
\otimes P_0,
\]
which is in $\mf I_3$. Again, since $\mf I_3$ is an ideal in $\mf T\otimes \mf 
T$, we see that both $P_0\otimes I$ and $I\otimes P_0$ are in $\mf I_3$ and 
hence $\mf I_3$ contains both $P_0\otimes \mf T$ and $\mf T\otimes P_0$. 
Therefore, 
$\mf K \otimes  \mf T +  \mf T \otimes \mf K$ is contained in $\mf I_3$. Now we 
consider 
the exact sequence
\[
\mf K \otimes  \mf T \oplus \mf T \otimes \mf K \lngarr  \mf T \otimes  \mf T 
\overset{\ds\Sigma}{\hbox to 25pt{\rightarrowfill}} C((\partial \bb D)^2) 
\lngarr 0
\]
and note that $\mf I_3$ is determined by $\Sigma(\mf I_3)$. Hence, we have $\mf 
I_3 = \{X\in \mf T\otimes  \mf T\mid \Sigma(X)|_{\mf Z_3} \equiv 0\}$ which 
implies 
that
\[
 \mf T(H^2(\bb D^2))/{\mf I_3} \cong C(\{(z_1,z_2)\in (\partial\bb D)^2 \mid 
z_1=z_2\}) \cong C(\partial\bb D).
\]
Therefore, there is no $C^*$-map from $\mf T(H^2(\bb D^2))$ to $\mf T(H^2(\bb 
D))$ that 
extends the correspondence $T_\varphi\to T_{\varphi|\mf Z_3}$ defined for 
$\varphi$ 
in $A(\bb D^2)$.

\section{Ideals in Toeplitz $\pmb{C^*}$-Algebras---General Case}\label{sec3}

\indent

The foregoing results might seem, at first, to be a little peculiar. For both 
of the cases in which the quotient Hilbert modules are well-behaved, namely 
examples one and three, there is no $C^*$-map, while in the second example there 
is a $C^*$-map but no quotient module map. As we said at the beginning of this 
note, at this point in our study we saw that we were attempting to answer the 
wrong question. However, there is a characterization of ideals in $\mf T(H^2(\bb 
D^2))$ which makes sense of the results about the foregoing examples.

We begin by analyzing the ideals in $\mf T \otimes\mf K$, one of the 
$C^*$-algebras 
in the diagram.

\begin{lem}\label{lem3.1}
The non-zero ideals $\mf I$ in $\mf T\otimes\mf K$ are parametrized by the 
closed 
subsets  $M$ of $\partial \bb D$ such that $\mf I_M = \{X\in \mf T\otimes\mf 
K\mid 
(\sigma\otimes 1)(X)(z) = 0$ for $z$ in $M\}$, where $\sigma\otimes 
1\colon \  \mf T \otimes\mf K \lngarr C(\partial\bb D) \otimes \mf K \cong 
C(\partial \bb D, \mf K)$.
\end{lem}

\begin{proof}
If $\mf I\ne (0)$, then $\mf K(H^2(\bb D^2))\cong \mf K \otimes \mf K \subset 
\mf I \subset \mf T \otimes \mf K$ and the result follows from the exactness of 
the 
sequence
\[
0 \lngarr \mf K\otimes \mf K \lngarr \mf T\otimes \mf K \overset{\ds 
\sigma\otimes 1}{\hbox to 25pt{\rightarrowfill}} C(\partial\bb D)\otimes\mf K 
\lngarr 0 
\]
and the structure of the ideals in $C(\partial\bb D) \otimes \mf K\cong 
C(\partial \bb D,\mf K)$.
\end{proof}

Next we consider the ideals in $C(\partial\bb D)\otimes \mf T$.

\begin{lem}\label{lem3.2}
A non-zero ideal $\mf J$ in $C(\partial\bb D) \otimes \mf T\cong C(\partial\bb 
D,\mf T)$ 
is characterized by two closed sets $Z' \subseteq \partial\bb D$ and $Z 
\subseteq (\partial\bb D)^2$ with the properties that
\begin{itemize}
\item[(1)] $Z'\times \partial\bb D \subset  Z$, and
\item[(2)] $\mf J = \{X\in C(\partial\bb D, \mf T) \mid X(z_2) = 0$ for $z_2$ is 
$Z'$ and
\[
\sigma[X(z_2)] (z_1) = 0\quad \text{for}\quad (z_1,z_2) \text{ in }  Z\}.
\]
\end{itemize}
\end{lem}

\begin{proof}
The set of values at $z_1$ in $\partial \bb D$ of the elements $X$ in $\mf J$ 
forms an 
ideal  in $\mf T$. Since the ideals in $\mf T$ are (0), $\mf K$, and the 
operators 
$\{X\in \mf T|\sigma(X)|_S \equiv 0\}$ for some closed subset $S$ of 
$\partial\bb D$, the result follows. The set $Z'$ corresponds to the subset 
of $\partial\bb D$ on which evaluation yields the zero ideal.
\end{proof}

An ideal in $\mf T\otimes \mf T$ gives rise to two ideals to which the preceding 
lemma
applies.

\begin{lem}\label{lem3.3}
If $\mf I$ is an ideal in $\mf T\otimes \mf T$, then the pairs of sets 
$(Z'_1,Z_1)$ and
$(Z'_2,Z_2)$ that characterize $(\sigma\otimes 1)(\mf I)$ and 
$(1\otimes\sigma)(\mf I)$, respectively, satisfy $Z_1 =  Z_2$.
\end{lem}

\begin{proof}
The result follows by identifying $Z_1$ and $Z_2$ with the zero set of 
the ideal $\Sigma(\mf I) \subset C((\partial\bb D)^2)$.
\end{proof}

We now build on this result to further analyze ideals in $\mf T\otimes \mf T$.

\begin{lem}\label{lem3.4}
If $\mf I$ is a non-zero ideal in $\mf T\otimes \mf T$, $\mf I_1$ is the ideal 
$\mf 
I\cap (\mf T\otimes \mf K)$ in $\mf T\otimes \mf K$ characterized by $M\subset 
\partial 
\bb D$, and $\mf I'$ is the ideal $(\sigma\otimes 1)(\mf I)$ in $C(\partial\bb 
D)\otimes \mf K$ characterized by the pair $(Z',Z)$ in $\partial\bb D$ 
and $(\partial\bb D)^2$, respectively, then $Z' = M$.
\end{lem}

\begin{proof}
If we let $\mf I'_1 = (\sigma\otimes 1)(\mf I_1)$, then $X$ is in $\mf I_1$ if 
and only if $(\sigma\otimes 1)(X)(z) = 0$ for $z$ in $M$. Thus, $\mf I'_1$ 
consists of the functions in $C(\partial\bb D)\otimes\mf K \cong C(\partial\bb 
D,\mf K)$ that vanish on $M$. Moreover, the set of values of the functions in 
$\mf I'_1$ at $z$ in $\partial\bb D$ forms an ideal $\mf I'_1(z)$ in the 
subalgebra $\mf I'(z)$ of values at $z$ of the functions in $\mf I'\subset 
C(\partial\bb D, \mf T)$. We claim that $\mf I'_1(z_0) = (0)$ for $z_0$ in 
$\partial\bb D$ implies $\mf I'(z_0) = (0)$ from which the result follows.

Suppose $X$ is an element of $\mf I$ for which $(\sigma\otimes 1)(X) (z_0)\ne 
0$. Since $\mf I$ is an ideal, by multiplying it by an appropriate element in 
$\mf T\otimes \mf T$ we can assume $(\sigma\otimes 1)(X)(z_0)$ is compact and 
nonzero. In 
particular, there exists a compact operator $K$ such that $(\sigma\otimes 
1)(X)(z_0)K\ne 0$. Then $X'=X(I\otimes  K)$ is in $\mf I$ as well as in $\mf 
T\otimes \mf K$ and thus in $\mf I'$. But $(\sigma\otimes 1)(X')(z_0) = 
(\sigma\otimes 1)(X)(z_0) K$ is a non-zero element of $\mf I'_1(z_0)$ which 
is a 
contradiction.
\end{proof}

We can now state the main result describing the ideals in $\mf T(H^2(\bb D^2))$.

\begin{thmb}
A non-zero ideal $\mf I$ in $\mf T(H^2(\bb D^2))$ is characterized by a triple 
$(Z', Z'',  Z)$ with $Z'$ and $Z''$ closed subsets of $\partial\bb 
D$ 
and $Z$ a closed subset of $(\partial\bb D)^2$ which contains both $Z'\times 
\partial\bb D$ and $\partial\bb D \times Z''$. Moreover,
\begin{gather*}
\mf I = \{X\in \mf T(H^2(\bb D^2))\mid \Sigma(X)|_{\mf Z} \equiv 0, 
(\sigma\otimes 
1)(X)(z) 
= 
0\\
\text{for $z$ is $Z'$ and } (1\otimes\sigma)(X)(z) = 0 \text{ for }
z \text{ in }  Z''\}.
\end{gather*}
Further, there is an isometry $\alpha = (\oplus \alpha_z) \oplus (\oplus\beta_z) 
\oplus 
\mu$ from $\mf T\otimes \mf T/\mf I$ to $\left(\bigoplus\limits_{z\in  Z'} 
\mf T\right) 
\oplus \left(\bigoplus\limits_{z\in  Z''} \mf T\right) \oplus C(Z)$, 
where 
$\alpha_z(X) = (\sigma\otimes 1)(X)(z)$, $\beta_z(X) = (1\otimes\sigma)(X)(z)$ 
and $\mu(X) = \Sigma(X)$. Finally, the range of $\alpha$ consists of the 
$\left(\bigoplus\limits_{z\in Z'} X_z\right) \oplus  
\left(\bigoplus\limits_{z\in  Z''} Y_z\right) \oplus f$, where $\{X_z\}$, 
$\{Y_z\}\subset \mf T$ and $f$ in $C((\partial\bb D)^2$ and 1)~~the 
maps 
$z\to X_z$ and $z\to Y_z$ are continuous from $Z'$ and $Z''$ to $\mf T$, 
respectively, and 2)~~$((1\otimes\sigma) X_z)(w) = (\sigma\otimes 1)Y_w)(z) = 
f(z,w)$ for $(z,w)$ in $Z'\times Z''$.
\end{thmb}

\begin{proof}
Using the preceding lemmas, it is possible to define the sets $Z', Z''$ 
and $Z$ as indicated. Namely, $Z$ is the zero set of the ideal 
$\Sigma(\mf I)$ in $C((\partial\bb D)^2)$, and $Z'$ and $Z''$ are the 
subsets of those $z$ in $\partial\bb D$ for which $(\sigma\otimes 1)(\mf I)(z) = 
(0)$ and $(1\otimes\sigma)(\mf I)(z) = (0)$, respectively. Clearly then $\mf I$ 
is contained in the ideal defined in the statement of the Theorem and our task 
is to show that they are equal. Hence, let $X$ be an operator in $\mf T\otimes 
\mf T$ 
that satisfies (1)~~$\Sigma(X)|_{\mf Z} \equiv 0$, (2)~~$(\sigma\otimes 1)(X)(z) 
= 0$ 
for $z$ in $Z'$ and (3)~~$(1\otimes\sigma)(X)(w) = 0$ for $w$ in $Z''$. We must 
show that $X$ is in $\mf I$.

Since $\Sigma(\mf I) = \{f\in C((\partial\bb D)^2)\mid f|_{\mf Z} \equiv 0\}$, 
there 
exists $X_0$ in $\mf I$ such that $\Sigma(X_0) = \Sigma(X)$ and hence 
$\Sigma(X-X_0) = 0$. Therefore, by Proposition \ref{pro2.2}, there exist $Y_1$ 
and $Y_2$ in $\mf T\otimes \mf K$ and $\mf K\otimes \mf T$, respectively, such 
that 
$X-X_0 = Y_1+Y_2$. If we evaluate the equation $(\sigma\otimes 1)(X-X_0) = 
(\sigma\otimes 1)(Y_1) + (\sigma\otimes 1)(Y_2)$ at $z$ on $Z'$, then we 
obtain $0 = (\sigma\otimes 1)(Y_1)(z) + (\sigma\otimes 1)(Y_2)(z)$ by the 
assumptions on $X$ and the definition of $Z'$. Since $(\sigma\otimes 1)(Y_2) 
= 0$, it follows that $(\sigma\otimes 1)(Y_1)(z) = 0$ for $z$ in $Z'$. 
Similarly we have $(1\otimes \sigma)(Y_1) =  0$ and $(1\otimes\sigma)(Y_2)(z) = 
0$ for $z$ in $Z''$. If we can show that $Y_1$ and $Y_2$ are in $\mf I$, 
then it will follow that $X$ is in $\mf I$ thereby completing the proof. By 
symmetry, it is enough to show that $Y_1$ is in $\mf I$.

By Lemma \ref{lem3.1}, it follows that $\mf I\cap (\mf T\otimes\mf K)$ is 
characterized by a closed subset $Z'_0$ of $\partial\bb D$ consisting of 
those $z$ for which $(\sigma\otimes 1)(\mf I)(z) = (0)$. Clearly, $Z' 
\subseteq  Z'_0$ and it is enough to show equality of these sets. For $z_0$ 
in $Z'_0\backslash Z'$ there exists an operator $W$ in $\mf I$ such that 
$(\sigma\otimes 1)(W)(z_0)\ne 0$. As in the proof of Lemma \ref{lem3.4}, one can 
show that one can obtain $W'$ in $\mf T\otimes\mf K$ for which $(\sigma\otimes 
1)(W')(z_0)\ne 0$ which completes the proof of the characterization of the 
ideals in $\mf T(H^2(\bb D^2))$.

To establish that $\alpha$ is an isometry we first we define $\bar\alpha$ as the 
direct 
sum of the maps $\bar\alpha_z(X) = (\sigma\otimes 1)(X)(z)$ for $z$ in 
$\partial\bb D$, $\bar\beta_w(X) = (1\otimes\sigma)(X)(w)$ for $w$ in 
$\partial\bb D$, and $\bar\mu(X) = \Sigma(X)$ for $X$ in $\mf T\otimes \mf T$. 
The range 
of both $\bar\alpha_z$ and $\bar\beta_w$ is $\mf T$, while the range of 
$\bar\mu$ 
is $C((\partial\bb D)^2)$. It remains to observe that 
$\bigcap\limits_{z\in\partial\bb D} \text{ ker } \bar\alpha_z = 
\text{ker}(\sigma\otimes 1) = \mf K\otimes \mf T$ and 
$\bigcap\limits_{w\in\partial\bb 
D} \text{ ker } \bar\beta_w = \text{ker}(1\otimes\sigma) = \mf T\otimes\mf K$.
Thus $\text{ker } \bar\alpha = (\mf T\otimes\mf K)\cap (\mf K\otimes \mf T) = 
\mf K 
\otimes \mf K$. If $\mf I$ is a non-zero ideal in $\mf T(H^2(\bb D^2))$, it 
contains 
$\mf K\otimes \mf K$.

Now consider $\overset{\ds=}{\alpha}\colon \ \mf T\otimes \mf T\lngarr 
\left(\bigoplus\limits_{z\in Z'} \mf T\right) \oplus 
\left(\bigoplus\limits_{w\in Z''}\mf T\right) \oplus C(Z)$, where $\mu(X) 
= 
\Sigma(X)|_{Z}$. Clearly the kernel of $\overset{\ds=}{\alpha}$ is $\mf I$ 
by the 
first statement of the Theorem and hence we can define $\alpha = 
\overset{\ds=}{\alpha} \circ\pi^{-1}$ on $\mf T\otimes \mf T/\mf I$, where $\pi$ 
is the 
quotient map from $\mf T\otimes \mf T$ to $\mf T\otimes \mf T/\mf I$. Since 
$\ker 
\overset{\ds=}{\alpha} = \mf I, \alpha$ is one-to-one and hence isometric.

Finally, we must show that the range of $\alpha$ consists of all elements $x 
= 
\left(\bigoplus\limits_{z\in Z'}X_z\right)\oplus 
\left(\bigoplus\limits_{w\in Z''}Y_w\right)\oplus f$ in 
$\left(\bigoplus\limits_{z\in Z'} \mf T\right)\oplus 
\left(\bigoplus\limits_{w\in Z''} \mf T\right)\oplus C(Z)$ that 
satisfying 
conditions (1) and (2). That all elements $x$ in the range of $\alpha$ 
satisfy
conditions (1) and (2) follows from properties of the maps in the tensor product 
diagram in Proposition \ref{pro2.1}. Hence, the remainder of the proof requires 
showing that for an element $x = \left(\bigoplus\limits_{z\in Z'} 
X_z\right) 
\oplus \left(\bigoplus\limits_{w\in Z''}Y_w\right)\oplus f$ satisfying (1) 
and (2), there exists an $X$ in $\mf T\otimes \mf T$ such that $\alpha_z(X) = 
X_z$ for 
$z$ in $Z'$, $\beta_w(X) = Y_w$ for $w$ in $Z''$ and 
$\Sigma(X) = f$.

Taking $\overline X$ in $\mf T\otimes \mf T$ such that $\Sigma(\overline X) = 
f$, the 
element 
\[
x' = \left(\bigoplus\limits_{z\in  Z'} X_z - \alpha_z(\overline 
X)\right) \oplus \left(\bigoplus\limits_{w\in Z''} Y_w - \beta_w(\overline 
X)\right) \oplus (f-\Sigma(\overline X))
\]
 also satisfies (1) and (2). Moreover, 
showing that $x'$ is in the range of $\alpha$ implies that $x$  is also. 
Hence we 
can assume for the given $x$ that $f\equiv 0$.

Assume we have such an element $x$. Extend the function $z\to X_z$ on $Z'$ to 
a function $\Phi$ defined from $\partial\bb D$ to $\mf T$ by extending it 
linearly 
on the open intervals in $\partial\bb D\backslash Z'$. Observe that $\Phi$ 
is continuous by (1) and hence in $C(\partial\bb D, \mf T) \cong C(\partial\bb 
D)\otimes \mf T$. Moreover, for $z'$ in $Z'$ and $w$ in $\partial\bb D$, we 
have
\[
((1\otimes\sigma)\Phi)(z,w) = ((1\otimes\sigma)X_z)(w) = f(z,w) = 0.
\]
Hence, $(1\otimes\sigma)X_z = 0$ for $z$ in $Z'$ and since $\Phi$ was 
defined 
by linear extension we have $(1\otimes\sigma) \Phi(z) = 0$ or 
$(1\otimes\sigma)\Phi = 0$. By the exactness of the sequence
\[
0 \lngarr C(\partial\bb D,\mf K) \lngarr C(\partial\bb D, \mf T)\lngarr 
C((\partial\bb D)^2)\lngarr 0,
\]
we see that $\Phi$ is in $C(\partial\bb D,\mf K)$ or $\Phi(z)$ is compact for 
$z$ in $\partial \bb D$. Since the sequence
\[
0\lngarr \mf K\otimes\mf K\lngarr \mf T\otimes \mf K \lngarr C(\partial\bb D) 
\otimes \mf K \lngarr 0
\]
is exact, there exists $\overset{\ds =}{X}$ in $\mf T\otimes \mf K$, and hence 
also 
in $\mf T\otimes \mf T$, such that $(\sigma\otimes 1)(\overset{\ds =}{X}) = 
\Phi$. 
Replace $x$ by $x-(\overset{\ds =}{X})$ to obtain an element $x'' = 
\left(\bigoplus\limits_{z\in Z'} X'_z\right) \oplus 
\left(\bigoplus\limits_{w\in Z''} Y'_w\right) \oplus f'$ in 
$\left(\bigoplus\limits_{z\in Z'} \mf T\right)\oplus 
\left(\bigoplus\limits_{w\in Z''} \mf T\right)\oplus C(Z)$ that satisfies 
(1) 
and (2) and such that $f'\equiv 0$ and $X'_z = X_z -(\sigma\otimes 
1)(\overset{\ds =}{X})(z) = X_z - \Phi(z) = 0$. Repeating the previous 
construction but reversing the roles of the first and second factors, we can 
replace $x''$ by an element $x''' = \left(\bigoplus\limits_{z\in  Z'} 
0\right) 
\oplus \left(\bigoplus\limits_{w\in Z''} 0\right)\oplus 0$ such that $x''' = 
x'' - \alpha (\overset{\ds\equiv}{X})$. Hence, $x''$ is in the range of 
$\alpha$ 
if and only if the zero element is, which concludes the proof.
\end{proof}

Examining the proof, one sees that the key is the fact that the intersection of 
$\ker(\sigma\otimes 1)$ and $\ker(1\otimes\sigma)$ in $\mf T\otimes \mf T$ is 
the 
intersection of $\mf K\otimes \mf T$ and $\mf T\otimes \mf K$ or $\mf K\otimes 
\mf K$, 
and that the latter ideal is contained in every nontrivial ideal in $\mf 
T\otimes \mf 
T$. If 
we consider the Toeplitz $C^*$-algebra over $H^2(\bb D^3)$, then one has a 
symbol map $\Sigma$ from $\mf T\otimes \mf T\otimes \mf T$ to $C((\partial \bb 
D)^3)$ and 
three symbol maps $\Sigma_1,\Sigma_2$, and $\Sigma_3$ from $\mf T\otimes \mf 
T\otimes \mf T$ 
to $C((\partial\bb D)^2, \mf T)$ with the intersection of their kernels again 
being 
$\mf K\otimes\mf K\otimes \mf K$. Thus one can formulate a description of the 
ideals in $\mf T(H^2((\partial\bb D)^3))$ analogous to the theorem in terms of 
$\Sigma,\Sigma_1$, $\Sigma_2$ and $\Sigma_3$. A similar approach will work for 
Toeplitz 
$C^*$-algebras defined for domains of rank greater than one. To some extent, 
this idea structure also shows up in index theory (cf.\ \cite{upmeier}).

\end{document}